\documentclass[12pt,reqno]{amsart}

\usepackage{amssymb,latexsym}

\usepackage{enumerate}
\allowdisplaybreaks
\usepackage[french,english]{babel}
\usepackage{amsmath}
\usepackage{graphicx}
\usepackage{amssymb}
\usepackage{bbm}
\usepackage{amsthm,mathtools}
\usepackage{ulem}
\usepackage{geometry}
\usepackage{tikz-cd}
\usepackage{mathrsfs}
\usepackage[colorinlistoftodos]{todonotes}
\usepackage{enumitem}
\usepackage{verbatim}
\usepackage[foot]{amsaddr}
\usepackage{dsfont}
\usepackage{cite}
\usepackage[T1]{fontenc}

\makeatletter

\@namedef{subjclassname@2010}{
	
	\textup{2020} Mathematics Subject Classification}

\makeatother
\newtheorem{thm}{Theorem}[section]
\newtheorem*{thm*}{Theorem}

\newtheorem{cor}[thm]{Corollary}

\newtheorem{lem}[thm]{Lemma}
\newtheorem{pro}[thm]{Proposition}
\theoremstyle{definition}

\numberwithin{equation}{section}

\newcommand{\bg}{\big}

\newcommand{\bgg}{\bigg}

\newcommand{\Bg}{\Big}

\newcommand{\Bgg}{\Bigg}

\newcommand{\lo}{\log_2}
\newcommand{\lt}{\log_3}

\newcommand{\inv}{^{-1}}

\newcommand{\mbr}{\mathbb{R}}

\newcommand{\mcp}{\mathcal{P}}

\newcommand{\mcl}{\mathcal{L}}

\newcommand{\mcz}{\mathcal{Z}}

\newcommand{\mcm}{\mathcal{M}}

\newcommand{\mmd}{\mathrm{d}}
\newcommand{\mme}{\mathrm{e}}
\newcommand{\mmi}{\mathrm{i}}
\newcommand{\mit}{\mathrm{i}t}
\newcommand{\mdt}{\mathrm{d}t}
\newcommand{\re}{\operatorname{Re}}
\newcommand{\im}{\operatorname{Im}}

\usepackage{hyperref}
\hypersetup{colorlinks=true,linkcolor=blue,anchorcolor=blue,citecolor=blue}
\usepackage{color}
\newcommand{\newabstract}[1]{%
	\par\bigskip
	\csname otherlanguage*\endcsname{#1}%
	\csname captions#1\endcsname
	\item[\hskip\labelsep\scshape\abstractname.]
}

\begin{document}

	\baselineskip=17pt

	\title[S. Zhao]{Two-Sided large oscillations of the argument of the Riemann zeta function and its iterates}

    \author{Shengbo Zhao\textsuperscript{1}}

	\address{1. School of Mathematical Sciences, Key Laboratory of Intelligent Computing and Applications (Tongji University), Ministry of Education, Tongji University, Shanghai 200092, China}

	\email{shengbozhao@hotmail.com}

	\begin{abstract} 
	   Assuming the Riemann hypothesis, we establish lower bounds for large oscillations of the iterated argument functions \(S_n(t)\) of the Riemann zeta function. As a corollary, we prove that, for every fixed integer \(n \geq 0\), 
       \[S_n(t)=\Omega_{\pm}\Bg(\frac{(\log t)^{1/2}(\lt t)^{1/2}}{(\lo t)^{n+1/2}}
\Bg).\]
For even \(n\), this recovers two results of Chirre and Mahatab in 2021. For odd \(n\), our result strengthens theirs by establishing the missing direction of oscillation with the same order of magnitude.
	\end{abstract}
	
	
	
	\maketitle

\section{Introduction}
\label{sec-int}

Throughout this paper, we denote the Riemann zeta function by \(\zeta(s)\), where \(s=\sigma+\mit\) is a complex variable. For \(1/2\le \sigma\le1\) and \(t >0\), define
\[S(\sigma,t):=\frac{1}{\pi}\arg\zeta(\sigma+\mit),\]
where the argument is obtained by continuous variation along straight line segments joining the points \(2\), \(2+\mit\) and \(\sigma+\mit\), assuming that the segment from \(\sigma+\mit\) to \(2+\mit\) has no zeros of \(\zeta(s)\), and with the convention that \(\arg\zeta(2)=0\). If this path has zeros of \(\zeta(s)\) (including the endpoint \(\sigma+\mit\)), define
\[S(\sigma,t)=\frac{S(\sigma,t+0)+S(\sigma,t-0)}{2}.\]
Set \(S_0(\sigma,t)=S(\sigma,t)\). Then, for integer \(n \ge 1\), define
\[S_n(\sigma,t)=\int_0^t S_{n-1}(\sigma,\tau)\mmd \tau + \delta_{n,\sigma},\]
where the constant \(\delta_{n,\sigma}\) is given by Fujii \cite{Fujii2002Comment}. More precisely, for \(n=2k\) with \(k\ge 1\),
\[
\delta_{2k,\sigma}=(-1)^{k-1}\int_\sigma^1 \int_{u_{2k}}^1 \cdots \int_{u_3}^1\int_{u_2}^1\mmd u_1\mmd u_2\cdots\mmd u_{2k} = \frac{(-1)^{k-1}(1-\sigma)^{2k}}{(2k)!},
\]
and for \(n=2k-1\) with \(k \ge 1\),
\[
\delta_{2k-1,\sigma}=\frac{(-1)^{k-1}}{\pi}\int_\sigma^\infty \int_{u_{2k-1}}^\infty \cdots \int_{u_3}^\infty\int_{u_2}^\infty\log|\zeta(u_1)|\mmd u_1\mmd u_2\cdots\mmd u_{2k-1}.
\]

The function \(S(t)\) is closely related to the distribution of the non-trivial zeros of the Riemann zeta function. The classical Riemann–von Mangoldt formula states that
\[N(t)=\frac{t}{2\pi}\log\Bg(\frac{t}{2\pi\mme}\Bg)+\frac{7}{8} +S(t) +O\bg(t\inv\bg),\]
where \(N(t)\) denotes the number of zeros \(\rho=\beta+\mmi\gamma\) of \(\zeta(s)\) with \(0<\gamma\leq t\), counted with multiplicity, where zeros with ordinate \(t\) are counted with half weight. In particular, for \(0<h \leq1 \), we have
\[N(t+h)-N(t)=\frac{h}{2\pi}\log t+S(t+h)-S(t)+O\bg(h+t\inv\bg).\]
This relation shows that the oscillations of \(S(t)\) are intimately connected with the local distribution of zeros of \(\zeta(s)\). More generally, the iterated functions \(S_n(t)\) describe successive integrals of these fluctuations. Understanding their large values and oscillations is therefore a natural problem in the study of the Riemann zeta function.

The study of large values of \(S(t)\) has a long history. In 1977, Montgomery \cite{Montgomery1977CommentMH} established, assuming the Riemann hypothesis (RH), the following \(\Omega_{\pm}\)-result:
\[S(t)=\Omega_{\pm}\Bg(\frac{(\log t)^{1/2}}{(\lo t)^{1/2}}\Bg).\]
Here and throughout this paper, we write \(\log_j\) for the \(j\)-th iterated logarithm. A heuristic prediction for the true order of magnitude of \(S(t)\) was proposed by Farmer, Gonek and Hughes \cite{Farmer2007Crelle}. Based on random matrix theory, they conjectured that
\[\limsup_{t\to\infty}\frac{S(t)}{(\log t)^{1/2}(\lo t)^{1/2}} = \frac{1}{\pi\sqrt{2}}.\]
This conjecture suggests that large values of \(S(t)\) may be substantially larger than the currently known lower bounds.

Tsang \cite{Tsang1986AA} subsequently obtained several refinements of the classical \(\Omega\)-results. More precisely, he proved that
\[S(t)=\Omega_\pm \Bg(\frac{(\log t)^{1/3}}{(\lo t)^{1/3}}\Bg),\quad \text{unconditionally},\quad S_1(t)=\Omega_\pm \Bg(\frac{(\log t)^{1/2}}{(\lo t)^{3/2}}\Bg),\quad\text{under RH}.\]
In addition, Tsang obtained unconditional estimates for \(S_1(t)\), with different orders in the two directions. In a subsequent paper, Tsang \cite{Tsang1991QJMath} further improved \(\Omega_+\)-result for \(S_1(t)\).

A significant development came from the resonance method. In 2018, Bondarenko and Seip \cite[Theorem 2]{Bondarenko2018MathAnn} proved, under RH, that
\[S(t)=\Omega\Bg(\frac{(\log t)^{1/2}(\lt t)^{1/2}}{(\lo t)^{1/2}}\Bg),\quad S_1(t)=\Omega_+\Bg(\frac{(\log t)^{1/2}(\lt t)^{1/2}}{(\lo t)^{3/2}}\Bg).\]
These results introduced an additional factor of \((\lt t)^{1/2}\) compared with the earlier conditional work, although the bound for \(S(t)\) did not distinguish between its positive and negative values.

Using a modified version of the resonance method in \cite{Bondarenko2018MathAnn}, Chirre \cite{Chirre2019JNT} investigated large values of \(S_n(\sigma,t)\) near the critical line. Assuming RH, when \(\sigma \in [1/2,1/2+1/(\lo T)]\), he obtained lower bounds of order, for sufficiently large \(T\),
\[\frac{(\log T)^{1-\sigma}(\lt T)^\sigma}{(\lo T)^{\sigma+n}}.\]
 When \(n\equiv 1 \pmod 4\), he established an \(\Omega_+\)-result for \(S_n(\sigma,t)\), whereas in the other cases, a corresponding \(\Omega\)-result is available. In 2021, Chirre and Mahatab subsequently made two important contributions. Specifically, in \cite{Chirre2021BLMS} they obtained that
 \[S(t)=\Omega_\pm\Bg(\frac{(\log t)^{1/2}(\lt t)^{1/2}}{(\lo t)^{1/2}}\Bg).\]
This result complements the missing negative direction in \cite{Bondarenko2018MathAnn}. On the other hand, Chirre and Mahatab \cite{Chirre2021JNT} extended this approach to the iterated functions \(S_n(\sigma,t)\) near the critical line. Assuming RH, they established large oscillations of \(S_n(\sigma,t+h)-S_n(\sigma,t)\) for every fixed integer \(n \ge 1\). Their result gives both directions when \(n\) is odd, but only one prescribed direction when \(n\) is even. Consequently, their results for \(S_n(t)\) are two-sided for even \(n\), whereas for odd \(n\) only one direction is obtained at this order of magnitude. They also derived certain unconditional one-sided consequences for higher iterates.

More recently, Inoue, Kobayashi and Toma \cite{Inoue2026MPCPS} further studied explicit extreme values of \(S(t+h)-S(t)\) in short intervals under RH. Their results also have applications to large and small gaps between consecutive zeros of the Riemann zeta function. For more developments and details about the argument of \(\zeta(s)\) and its iterates, we recommend \cite{Bui2014JMAA,Carneiro2013MathAnn,Carneiro2018MPCMS,Carneiro2019PublM,Fujii2002Comment,Inoue2022AA} and the references therein.

Motivated by these developments, we investigate whether large oscillation estimates of Chirre and Mahatab \cite{Chirre2021JNT} can be strengthened to both directions for every fixed integer \(n\geq1\). Our main result gives an affirmative answer under RH, without imposing any restriction on the parity of \(n\).

\begin{thm}
    \label{thm1}
    Assume RH. Let \(n \ge 1\) be a fixed integer and let \(0<\beta<1\) be fixed. Then there exists a constant \(c_{n,\beta}>0\) such that, for all sufficiently large \(T\) and \(h \in \bg[0,(\lo T)\inv\bg]\), we have
    \begin{align*}
        \max_{T^\beta \le t \le T} \varepsilon \bg(S_n(t+h)-S_n(t)\bg) \ge c_{n,\beta} h\frac{(\log T)^{1/2}(\lt T)^{1/2}}{(\lo T)^{n-1/2}},\quad \varepsilon\in\{-1,1\}.
    \end{align*}
\end{thm}

\begin{cor}
    \label{cor2}
    Assume RH. Let \(n \ge 0\) be a fixed integer and let \(0<\beta<1\) be fixed. Then there exists a constant \(c_{n,\beta}>0\) such that, for all sufficiently large \(T\), we have
    \[
    \max_{T^\beta \le t \le T} \varepsilon S_n(t) \ge c_{n,\beta} \frac{(\log T)^{1/2}(\lt T)^{1/2}}{(\lo T)^{n+1/2}},\quad \varepsilon\in\{-1,1\}.
    \]
    In particular,
    \[S_n(t)=\Omega_{\pm}\Bg(\frac{(\log t)^{1/2}(\lt t)^{1/2}}{(\lo t)^{n+1/2}}
\Bg).\]
\end{cor}

We now compare our results with the work of Chirre and Mahatab. In fact, \cite[Theorem 1]{Chirre2021JNT} establishes two-sided large oscillations of \(S_n(t+h)-S_n(t)\) when \(n\) is odd, but only one prescribed direction when \(n\) is even. Theorem \ref{thm1} removes this restriction and establishes both directions for every fixed integer \(n\geq1\), with the same order of magnitude. On the other hand, for \(S_n(t)\), Corollary \ref{cor2} recovers \cite[Corollary 2]{Chirre2021JNT} when \(n\) is even and \(\sigma=1/2\). When \(n\) is odd, it strengthens their result by establishing the missing direction. Furthermore, when \(n=0\), Corollary \ref{cor2} recovers \cite[Theorem 1]{Chirre2021BLMS}.

Although our results are stated on the critical line for simplicity, the same method yields the corresponding estimates uniformly for \(1/2\leq\sigma\leq1/2+(\log_2T)^{-1}\), with \(S_n(t)\) replaced by \(S_n(\sigma,t)\).

The choice of the kernel \(K(u)\), defined in \eqref{linear}, is inspired by Jaskari \cite{Jaskari2024FACM}. The resulting convolution formula, established in Proposition \ref{prop}, plays a crucial role in the proof of Theorem \ref{thm1}. In contrast to \cite[Proposition 8]{Chirre2021JNT}, our approach avoids determining which of two finite sums yields the main term based on the parity of \(n\). This is achieved by suitably choosing the phase of the sine function in \(K(u)\). Since \(\mmi^n(-\mmi)^{n-1}=\mmi\), this choice guarantees that, after taking the imaginary part, the second sum always appears with coefficient \(\pm1/2\), independent of \(n\).

\section{Preliminary tools}
\label{sec-pre}

\subsection{A convolution formula}

For \(h\in\mbr\) and a function \(F\), we write
\[\Delta_h F(t)=F(t+h)-F(t-h).\]
The following two lemmas, due to Chirre and Mahatab \cite{Chirre2021JNT}, play a crucial role in establishing the convolution formula we require. 

\begin{lem}[Chirre-Mahatab \cite{Chirre2021JNT}, Lemma 6]
\label{Chirre-lem-6}
    Assume RH. Let \(0<\beta<1\) be a fixed number. Let \(\alpha>0\), \(H \in \mbr\) and \(T\) be sufficiently large. Then for \(1/2<v\le2\), \(0 \le h \le 1\) and \(T^\beta \le t \le T\log T\), we have
    \begin{align*}
    \int_{-(\log T)^3}^{(\log T)^3} &\Delta_h \log\zeta(v+\mmi(t+u))\Bg(\frac{\sin\alpha u}{u}\Bg)^2 \mme^{\mmi Hu}\mmd u \\
    &= - \pi\mmi \sum_{m \ge 2}\frac{\Lambda(m)\omega_m(\alpha,H)\sin(h\log m)}{(\log m)m^{v+\mit}} + O\Bg(\frac{h(1+|\log(2v-1)|)\mme^{2\alpha+|H|}}{(\log T)^3} \Bg),
    \end{align*}
    where \(\omega_m(\alpha,H)=\max\{0,2\alpha-|H-\log m|\}\).
\end{lem}

\begin{lem}[Chirre-Mahatab \cite{Chirre2021JNT}, Lemma 7]
\label{Chirre-lem-7}
    Assume RH. Let \(n \ge 1\) be a fixed integer, and \(t,h \in \mbr\) such that \(t \neq \pm h\). Then, we have
    \[
    \Delta_h S_n(t) = \frac{1}{\pi} \im \bgg( \frac{\mmi^n}{(n-1)!}\int_{1/2}^2\Bg(v-\frac{1}{2}\Bg)^{n-1}\Delta_h\log\zeta(v+\mit)\mmd v\bgg)+O_n(h).
    \]
\end{lem}

We now construct the convolution formula we need. To this end, for large \(T\), \(n \ge 1\) and \(\varepsilon \in \{-1,1\}\), define
\begin{align}
    \label{linear}
    K(u)&=2+\varepsilon\sin\Bg(\frac{n-1}{2}\pi-u\lo T\Bg)\nonumber\\
    &=2\mme^0+\frac{\varepsilon\mmi}{2}(-\mmi)^{n-1}\mme^{\mmi u \lo T}-\frac{\varepsilon\mmi^n}{2}\mme^{-\mmi u \lo T}.
\end{align}
Note that \(K(u)\) is always positive. More precisely, \(1\le K(u) \le 3\).
\begin{pro}
\label{prop}
    Assume RH. Let \(n \ge 1\) be a fixed integer and let \(0<\beta<1\) and \(0 < \gamma \le 1/2\) be fixed. Let \(T\) be sufficiently large. Then for \(0 \le h \le 1\) and \(T^\beta \le t \le T\log T\), we have
    \begin{align*}
        \int_{-(\log T)^3}^{(\log T)^3}\Delta_h &S_n(t+u)\Bg(\frac{\sin(\gamma u \lo T)}{u}\Bg)^2 K(u) \mmd u \\
        &= \im\bg(2\mmi^{n+3}A_1(n,t)\bg)+\frac{\varepsilon}{2}\re A_2(n,t) + O_n(h\lo T),
    \end{align*}
    where
    \[A_j(n,t)= \sum_{m \ge 2}\frac{\Lambda(m)a_j(m)}{(\log m)^{n+1}m^{1/2+\mit}}, \quad j=1,2.\]
    Moreover,
    \begin{align*}
        &a_1(m):=a_1(m,T,h) = \omega_m(\gamma\lo T,0)\sin(h\log m), \\
        &a_2(m):=a_2(m,T,h) = \omega_m(\gamma\lo T,\lo T)\sin(h\log m).
    \end{align*}
\end{pro}

\begin{proof}
    Set \(\alpha=\gamma \lo T\) in Lemma \ref{Chirre-lem-6}, and set \(H=-\lo T, 0, \lo T\), respectively. Then, we have
    \begin{align}
    \label{prop-1}
        \int_{-(\log T)^3}^{(\log T)^3}\Delta_h &\log\zeta(v+\mmi(t+u))\Bg(\frac{\sin(\gamma u \lo T)}{u}\Bg)^2 K(u) \mmd u \nonumber \\
        &= -2\pi\mmi \sum_{m \ge 2}\frac{\Lambda(m)a_1(m)}{(\log m)m^{v+\mit}} +\frac{\pi\varepsilon}{2}(-\mmi)^{n-1} \sum_{m \ge 2}\frac{\Lambda(m)a_2(m)}{(\log m)m^{v+\mit}} \nonumber \\
        &\quad+O\Bg(\frac{h(1+|\log(2v-1)|)}{(\log T)^{2-2\gamma}} \Bg).
    \end{align}
    Here, we use the linear combination \eqref{linear}. Furthermore, when \(0<\gamma\le 1/2\), we have \(\omega_m(\gamma\lo T,-\lo T)=0\) for all \(m \ge 2\). Multiplying both sides of \eqref{prop-1} by \((v-1/2)^{n-1}\) and integrating over \(v \in [1/2,2]\) gives 
    \begin{align}
        \label{prop-2}
        \int_{-(\log T)^3}^{(\log T)^3}&\bgg(\int_{1/2}^2\Bg(v-\frac{1}{2}\Bg)^{n-1}\Delta_h  \log\zeta(v+\mmi(t+u)) \mmd v\bgg)\Bg(\frac{\sin(\gamma u \lo T)}{u}\Bg)^2 K(u) \mmd u \nonumber \\
        &= -\sum_{m \ge 2}\frac{\Lambda(m)}{(\log m)m^{\mit}} \bgg(\Bg(2\pi\mmi a_1(m)-\frac{\pi\varepsilon}{2}(-\mmi)^{n-1}a_2(m)\Bg)\int_{1/2}^2 \frac{(v-1/2)^{n-1}}{m^v}\mmd v \bgg) \nonumber\\
        &\quad +O\bgg(\frac{h}{(\log T)^{2-2\gamma}}\int_{1/2}^2 \Bg(v-\frac{1}{2}\Bg)^{n-1}(1+|\log(2v-1)|\mmd v \bgg) \nonumber \\
        &= -\sum_{m \ge 2}\frac{\Lambda(m)}{(\log m)m^{\mit}} \bgg(\Bg(2\pi\mmi a_1(m)-\frac{\pi\varepsilon}{2}(-\mmi)^{n-1}a_2(m)\Bg)\int_{1/2}^2 \frac{(v-1/2)^{n-1}}{m^v}\mmd v \bgg) \nonumber\\
        &\quad +O\Bg(\frac{h}{(\log T)^{2-2\gamma}}\Bg).
    \end{align}
    From the argument presented in \cite[p. 250]{Chirre2021JNT}, the inner integral contained in the main term on the right-hand side of \eqref{prop-2} satisfies, for each \(m \ge 2\), 
    \[\int_{1/2}^2 \frac{(v-1/2)^{n-1}}{m^v}\mmd v = \frac{(n-1)!}{(\log m)^n\sqrt{m}}+O\Bg(\frac{1}{m^2\log m}\Bg) .\]
    Combining this estimate with \eqref{prop-2}, we get
    \begin{align*}
        \int_{-(\log T)^3}^{(\log T)^3}&\bgg(\int_{1/2}^2\Bg(v-\frac{1}{2}\Bg)^{n-1}\Delta_h  \log\zeta(v+\mmi(t+u)) \mmd v\bgg)\Bg(\frac{\sin(\gamma u \lo T)}{u}\Bg)^2 K(u) \mmd u \\
        &= -2\pi(n-1)! \mmi A_1(n,t)+ \frac{\pi\varepsilon}{2}(n-1)!(-\mmi)^{n-1} A_2(n,t)+ O_n(h\lo T).
    \end{align*}
    Here we also use the following trivial upper bound that, for \(j=1,2\), 
    \begin{align}
        \label{ajmupper}
        |a_j(m)| \ll h \log m \lo T. 
    \end{align}
    Multiplying by the factor \(\mmi^n/(\pi(n-1)!)\) and taking the imaginary part of both sides, we complete the proof of Proposition \ref{prop}.
\end{proof}

\subsection{The resonator and some estimates}

We now construct the required resonator. Indeed, the resonator we use here coincides with the one chosen in \cite{Chirre2021BLMS}, both following the idea of \cite{Bondarenko2017DukeMJ,Bondarenko2018MathAnn,Chirre2019JNT}. For completeness, we recall its construction and the estimates needed in our proof.

Fix \(\beta \in (0,1)\) and set \(\kappa=(1-\beta)/2\), so that \(\kappa+\beta<1\). For sufficiently large \(T\), set \(N=\lfloor T^\kappa \rfloor\). Let \(\mcp\) denote the set of primes \(p\) satisfying
\[\mme \log N \lo N<p \le \exp\bg((\lo N)^{1/8}\bg)\log N \lo N.\]
Next, we define the multiplicative function $f(n)$ to be supported on the set of square-free numbers, with values for $p \in \mathcal{P}$ as
\[f(p)=\sqrt{\frac{\log N \lo N}{\lt N}}\frac{1}{\sqrt{p}(\log p-\lo N-\lt N)}.\]
Furthermore, if $p \notin \mathcal{P}$, set $f(p)=0$. 

For \(k \in \{1,\dots,\lfloor(\lo N)^{1/8\rfloor}\}\), we define the set $\mathcal{P}_{k}$ of prime numbers as follows
\[\mcp_k :=\bg\{p : \mme^k \log N \log_{2}N < p \leq \mme^{k+1} \log N\log_{2}N \bg\}\]
and
\[\mcm_k :=\Bg\{n \in \mathrm{supp}(f): n \,\, \text{has at least}\,\, \frac{3\log N}{k^2\lt N}\,\, \text{prime divisors in}\,\, \mcp_k\Bg\}.\]
Furthermore, define
\[\mcm := \mathrm{supp}(f) \setminus \bigcup_{k = 1}^{\lfloor(\log_2 N)^{1/8}\rfloor} \mcm_{k}.\]
Let \(\mcl\) be the set of integers \(\ell\) such that 
\[
\Bg[\Bg(1+\frac{1}{T}\Bg)^\ell,\Bg(1+\frac{1}{T}\Bg)^{\ell + 1}\Bg] \bigcap \mcm \neq \emptyset.
\]
Then, we let \(m_\ell\) be the minimum of \(\bg[(1+T\inv)^\ell, (1+T\inv)^{\ell+ 1}\bg] \cap \mcm\) for \(\ell\in\mcl\). Define the set
\[\mcm^\prime := \{m_\ell:\ell \in \mcl \},\]
and define
\[r(m_\ell):=\Bgg(\sum_{\substack{(1-T\inv)^{\ell-1}\le n \le(1+T\inv)^{\ell+2}\\n \in \mcm}}f(n)^2\Bgg)^{1/2}\]
for every \(m_\ell \in \mcm^\prime\). Moreover, the resonator is defined by
\[R(t):= \sum_{m\in\mcm^\prime}r(m)m^{-\mit}.\]

We conclude this section with some estimates. The Cauchy-Schwarz inequality implies that
\begin{align}
    \label{Rupper}
    |R(t)|^2 \le R(0)^2 \ll T^\kappa\sum_{k \in \mcm}f(k)^2.
\end{align}
Furthermore, set \(\Phi(y):=\mme^{-y^2/2}\) as in \cite{Bondarenko2017DukeMJ}. It follows from \cite[Lemma 5]{Bondarenko2018MathAnn} (also see \cite[Proposition 4]{Chirre2021BLMS}) that
\begin{align}
    \label{momentupper}
    \int_{-\infty}^\infty|R(t)|^2\Phi\Bg(\frac{t}{T}\Bg)\mdt \ll T\sum_{k \in \mcm}f(k)^2.
\end{align}
We will make frequent use of \eqref{Rupper} and \eqref{momentupper} in the subsequent argument. Finally, owing to the work of Bondarenko and Seip \cite{Bondarenko2018MathAnn}, we state the following lemma:
\begin{lem}[Bondarenko-Seip \cite{Bondarenko2018MathAnn}, Lemma 7]
\label{Seip-lem-7}
    If
    \[G(t) = \sum_{m\ge 2}\frac{\Lambda(m)u_m}{(\log m)m^{1/2+\mit}}\]
    is absolutely convergent and \(u_m \ge 0\) for each \(m \ge 2\), then we have
    \[
    \int_{-\infty}^\infty G(t)|R(t)|^2\Phi\Bg(\frac{t}{T}\Bg)\mdt \gg T \frac{(\log T)^{1/2} (\lt T)^{1/2}}{(\lo T)^{1/2}} \Bg(\min_{p \in \mcp}u_p \Bg)\sum_{k\in \mcm}f(k)^2.
    \]
\end{lem}

\section{Proof of Theorem \ref{thm1}}
\label{sec-pthm1}

In this section, we combine the resonance method with Proposition \ref{prop} to prove Theorem \ref{thm1}. The case \(h=0\) is trivial, so we assume that \(h \in (0,(2\lo T)\inv]\) throughout the proof. Set \(\gamma=1/8\). First, by Proposition \ref{prop}, we have
\begin{align}
\label{long}
    &\int_{T^\beta}^{T\log T}\bgg(\int_{-(\log T)^3}^{(\log T)^3}\Delta_h S_n(t+u)\Bg(\frac{\sin(\gamma u \lo T)}{u}\Bg)^2 K(u) \mmd u\bgg)|R(t)|^2\Phi\Bg(\frac{t}{T}\Bg)\mdt \nonumber\\
    &\quad = 2\im\bgg(\mmi^{n+3}\sum_{m \ge 2}\frac{\Lambda(m)a_1(m)}{(\log m)^{n+1}\sqrt{m}}\Bg(\int_{T^\beta}^{T\log T} m^{-\mit}|R(t)|^2\Phi\Bg(\frac{t}{T}\Bg)\mdt\Bg)\bgg) \nonumber \\
    &\qquad + \frac{\varepsilon}{2}\re\bgg(\sum_{m \ge 2}\frac{\Lambda(m)a_2(m)}{(\log m)^{n+1}\sqrt{m}}\Bg(\int_{T^\beta}^{T\log T} m^{-\mit}|R(t)|^2\Phi\Bg(\frac{t}{T}\Bg)\mdt\Bg)\bgg) \nonumber\\
    &\qquad +O\Bg(hT\lo T\sum_{k \in \mcm}f(k)^2\Bg).
\end{align}
Here, in the last step, we make use of \eqref{momentupper} and the fact that both \(|R(t)|^2\) and \(\Phi(t/T)\) are non-negative.

We first consider the second term on the right-hand side of \eqref{long}, namely the one associated with \(a_2(m)\), as this contributes the main term. To this end, we wish to expand the integral over \(\mathbb R\), so that Lemma \ref{Seip-lem-7} can be applied. By the definition of \(\omega_m\), \(a_2(m)=0\) unless \((\log T)^{3/4} \le m \le (\log T)^{5/4}\). Thus, combining this with \eqref{ajmupper} and \eqref{Rupper}, we obtain
\begin{align}
\label{a2upper-1}
    \bgg|\sum_{m \ge 2}\frac{\Lambda(m)a_2(m)}{(\log m)^{n+1}\sqrt{m}}&\bgg(\int_{0}^{T^\beta} m^{-\mit}|R(t)|^2\Phi\Bg(\frac{t}{T}\Bg)\mdt\bgg)\bgg| \nonumber\\
    &\quad \ll hT^{\kappa+\beta}\lo T \sum_{k \in \mcm}f(k)^2 \sum_{2 \le m \le(\log T)^{5/4}} \frac{\Lambda(m)}{(\log m)^{n+1}\sqrt{m}}.
\end{align}
By using \cite[Lemma B.1]{Carneiro2019PublM}, we have
\[\sum_{2 \le m \le(\log T)^{5/4}} \frac{\Lambda(m)}{(\log m)^{n+1}\sqrt{m}} \ll \frac{(\log T)^{5/8}}{(\lo T)^n}.\]
Since \(\kappa+\beta<1\), substituting this upper bound into \eqref{a2upper-1} implies that 
\begin{align}
    \label{tail-1}
    &\bgg|\sum_{m \ge 2}\frac{\Lambda(m)a_2(m)}{(\log m)^{n+1}\sqrt{m}}\Bg(\int_{0}^{T^\beta} m^{-\mit}|R(t)|^2\Phi\Bg(\frac{t}{T}\Bg)\mdt\Bg)\bgg| \nonumber\\
    &\quad \ll hT^{\kappa+\beta}\frac{(\log T)^{5/8}}{(\lo T)^{n-1}} \sum_{k \in \mcm}f(k)^2 \ll hT\sum_{k \in \mcm}f(k)^2.
\end{align}
Moreover, the exponential decay of \(\Phi\) yields the following crude upper bound:
\begin{align}
    \label{tail-2}
    \bgg|\sum_{m \ge 2}\frac{\Lambda(m)a_2(m)}{(\log m)^{n+1}\sqrt{m}}\Bg(\int_{T\log T}^\infty m^{-\mit}|R(t)|^2\Phi\Bg(\frac{t}{T}\Bg)\mdt\Bg)\bgg| \ll  hT\sum_{k \in \mcm}f(k)^2.
\end{align}

On the other hand, we seek an upper bound for the first term on the right-hand side of \eqref{long}. By the definition of \(\omega_m\), we have \(a_1(m) = 0\) whenever \(m >(\log T)^{1/4}\). Hence, the sum runs over \(2 \le m \le (\log T)^{1/4}\). Together with \eqref{ajmupper} and \eqref{momentupper}, this gives
\begin{align}
    \label{A1upper-1}
    &\bgg|\sum_{m\ge2}\frac{\Lambda(m)a_1(m)}{(\log m)^{n+1}\sqrt{m}}\Bg(\int_{T^\beta}^{T\log T} m^{-\mit}|R(t)|^2\Phi\Bg(\frac{t}{T}\Bg)\mdt\Bg)\bgg| \nonumber\\
    &\quad \le\bgg|h\lo T\sum_{2\le m \le(\log T)^{1/4}}\frac{\Lambda(m)}{(\log m)^n \sqrt{m}} \int_{-\infty}^\infty|R(t)|^2\Phi\Bg(\frac{t}{T}\Bg)\mdt \bgg| \nonumber\\
    &\quad \ll hT\lo T\sum_{2\le m \le(\log T)^{1/4}}\frac{\Lambda(m)}{(\log m)^n \sqrt{m}}\sum_{k \in \mcm}f(k)^2. 
\end{align}
By using \cite[Lemma B.1]{Carneiro2019PublM} again, we have
\[\sum_{2\le m \le(\log T)^{1/4}}\frac{\Lambda(m)}{(\log m)^n \sqrt{m}} \ll \frac{(\log T)^{1/8}}{(\lo T)^n}.\]
Thus, inserting this into \eqref{A1upper-1} yields that
\begin{align}
    \label{A1upper}
    \bgg|\sum_{m\ge2}\frac{\Lambda(m)a_1(m)}{(\log m)^{n+1}\sqrt{m}}\Bg(\int_{T^\beta}^{T\log T} m^{-\mit}|R(t)|^2\Phi\Bg(\frac{t}{T}\Bg)\mdt\Bg)\bgg| \ll h \frac{T(\log T)^{1/8}}{(\lo T)^{n-1}}\sum_{k \in \mcm}f(k)^2. 
\end{align}

To simplify the notation, we write
\[\mcz=\Bg(\max_{T^\beta/2 \le t \le 2T\log T} \varepsilon\Delta_h S_n(t) \Bg),\quad \mcz_+=\max\{0,\mcz\}.\]
We shall prove that \(\mcz_+>0\) for sufficiently large \(T\), which will imply that \(\mcz=\mcz_+\). 
Since
\[
T^\beta-(\log T)^3\geq\frac{T^\beta}{2},\quad T\log T+(\log T)^3\leq 2 T\log T
\]
for sufficiently large \(T\), we have
\[
\varepsilon \Delta_h S_n(t+u)\leq \mcz_+
\]
whenever \(T^\beta\leq t\leq T\log T\) and \(|u|\leq(\log T)^3\). Using \(1\leq K(u)\le 3\), together with the classical identity
\[\int_{-\infty}^\infty \Bg(\frac{\sin(au)}{u}\Bg)^2\mmd u=\pi a,\]
we obtain
\begin{align}
    \label{LHS1}
    &\quad\varepsilon\int_{T^\beta}^{T\log T}\bgg(\int_{-(\log T)^3}^{(\log T)^3}\Delta_h S_n(t+u)\Bg(\frac{\sin(u\lo T/8)}{u}\Bg)^2 K(u)\mmd u\bgg)|R(t)|^2\Phi\Bg(\frac{t}{T}\Bg)\mdt \nonumber \\
    & \le \mcz_+\int_{-(\log T)^3}^{(\log T)^3}\Bg(\frac{\sin(u\lo T/8)}{u}\Bg)^2 K(u)\mmd u \int_{T^\beta}^{T\log T}|R(t)|^2\Phi\Bg(\frac{t}{T}\Bg)\mdt \nonumber \\
    & \le \frac{3\pi}{8}\mcz_+ \lo T \int_{-\infty}^\infty |R(t)|^2\Phi\Bg(\frac{t}{T}\Bg)\mdt \ll \mcz_+ T\lo T \sum_{k \in \mcm}f(k)^2.
\end{align}
Here, in the last step, we apply \eqref{momentupper}. Noting that both \(|R(t)|^2\) and \(\Phi(t/T)\) are real and even, we have
\begin{align}
    \label{even}
    \re\bgg(\int_0^\infty m^{-\mit}|R(t)|^2\Phi\Bg(\frac{t}{T}\Bg)\mdt\bgg)=\frac{1}{2}\int_{-\infty}^\infty m^{-\mit}|R(t)|^2\Phi\Bg(\frac{t}{T}\Bg)\mdt.
\end{align}
Multiplying \eqref{long} by \(\varepsilon\), and then using \eqref{tail-1}, \eqref{tail-2}, \eqref{A1upper}, \eqref{LHS1} and \eqref{even}, we obtain
\begin{align}
    \label{max1}
    \mcz_+ T \lo T \sum_{k \in \mcm}f(k)^2 &\gg \frac{1}{4}\sum_{m \ge 2}\frac{\Lambda(m)a_2(m)}{(\log m)^{n+1}\sqrt{m}}\int_{-\infty}^\infty m^{-\mit}|R(t)|^2\Phi\Bg(\frac{t}{T}\Bg)\mdt\nonumber  \\
    &\quad+ O\Bg(hT(\log T)^{1/8}\lo T\sum_{k \in \mcm}f(k)^2\Bg).
\end{align}
We next show that the main term on the right-hand side of \eqref{max1} is positive. Since we assume \(h \le (\lo T)\inv\), it follows that 
\[\sin(h\log m) \gg h\log m\]
whenever \(m \le (\log T)^2\). Moreover, by the definition of \(\omega_m\), 
\[\omega_m\Bg(\frac{1}{8}\lo T,\lo T\Bg) \gg \log m\]
holds for all \(m \in [(\log T)^{7/8},(\log T)^{9/8}]\). Hence, when \(m \in [(\log T)^{7/8},(\log T)^{9/8}]\),
\[a_2(m) \gg h(\log m)^2.\]
Since \(N=\lfloor T^\kappa \rfloor\), we deduce that
\[(\log T)^{7/8} \le p \le (\log T)^{9/8}\]
for each \(p \in \mcp\) when \(T\) is sufficiently large. Thus, we get \(a_2(p) \gg h(\log p)^2\) and 
\[\min_{p \in \mcp} \frac{a_2(p)}{(\log p)^n} \gg \min_{p \in \mcp} \frac{h}{(\log p)^{n-2}} \gg \frac{h}{(\lo T)^{n-2}}.\]
Then, setting \(u_m=a_2(m)/(\log m)^n\) in Lemma \ref{Seip-lem-7} yields 
\begin{align}
    \label{laterG}
    \sum_{m \ge 2}\frac{\Lambda(m)a_2(m)}{(\log m)^{n+1}\sqrt{m}}\int_{-\infty}^\infty m^{-\mit}&|R(t)|^2\Phi\Bg(\frac{t}{T}\Bg)\mdt \nonumber\\
    &\gg hT\frac{(\log T)^{1/2}(\lt T)^{1/2}}{(\lo T)^{n-3/2}}\sum_{k \in \mcm}f(k)^2.
\end{align}
Combining \eqref{max1} with \eqref{laterG}, we observe that
\[\frac{hT(\log T)^{1/8}\lo T}{hT(\log T)^{1/2}(\lt T)^{1/2}(\lo T)^{3/2-n}} = \frac{(\lo T)^{n-1/2}}{(\log T)^{3/8}(\lt T)^{1/2}}=o(1),\]
since \(n\) is fixed. Thus, for all sufficiently large \(T\), \eqref{max1} and \eqref{laterG} imply
\[\mcz_+ \gg_{n,\beta} h\frac{(\log T)^{1/2}(\lt T)^{1/2}}{(\lo T)^{n-1/2}}>0.\]
Consequently, \(\mcz_+=\mcz\), and hence
\[ \max_{T^\beta/2 \le t \le 2T\log T} \varepsilon\Delta_h S_n(t) \gg_{n,\beta} h \frac{(\log T)^{1/2}(\lt T)^{1/2}}{(\lo T)^{n-1/2}}.\]

Finally, we adjust the parameters to obtain the desired result. Recall that
\[
\Delta_{h/2}S_n(v) =S_n(v+h/2)-S_n(v-h/2).
\]
Putting \(t=v-h/2\), we have
\[
\Delta_{h/2}S_n(v)=S_n(t+h)-S_n(t).
\]
Now, fix \(0<\beta<\beta_1<1\) and put \(T_1 = T/(4\log T)\). We apply the preceding estimate with \(T\) replaced by \(T_1\), \(\beta\) replaced by \(\beta_1\), and \(h\) replaced by \(h/2\). Since \(0<h\leq(\log_2T)^{-1}\), we have
\[0 < \frac{h}{2} \le \frac{1}{2\lo T_1}\]
for sufficiently large \(T\). Thus,
\[
\max_{T_1^{\beta_1}/2\leq v\leq 2T_1\log T_1} \varepsilon\Delta_{h/2}S_n(v)
\gg_{n,\beta_1} h\frac{(\log T_1)^{1/2}(\log_3 T_1)^{1/2}}
{(\log_2 T_1)^{n-1/2}}.
\]
Moreover, since \(\beta_1>\beta\), we have, for all sufficiently large \(T\),
\[\Bg[\frac{T_1^{\beta_1}}{2}-\frac{h}{2},2T_1 \log T_1-\frac{h}{2}\Bg] \subseteq [T^\beta,T]\]
Consequently,
\begin{align*}
    \max_{T^\beta\leq t\leq T} \varepsilon\big(S_n(t+h)-S_n(t)\big)  \gg_{n,\beta}h\frac{(\log T)^{1/2}(\log_3 T)^{1/2}}{(\log_2 T)^{n-1/2}},
\end{align*}
where we used \(\log_j T_1\sim\log_j T\) for \(j=1,2,3\). Since the estimate holds for each \(\varepsilon\in\{-1,1\}\), the proof of Theorem \ref{thm1} is complete.


\section{Proof of Corollary \ref{cor2}}
\label{sec-pcor2}

Set \(\beta^\prime=(1+\beta)/2\). Applying Theorem \ref{thm1} with \(n\) replaced by \(n+1\), \(T\) replaced by \(T/2\), we choose \(h=(\lo (T/2))^{-1}\), which lies in the admissible range of the theorem. Then, we obtain
\[\max_{(T/2)^{\beta'}\leq t\leq T/2} \pm \big(S_{n+1}(t+h)-S_{n+1}(t)\big) \gg_{n,\beta}h \frac{(\log T)^{1/2}(\lt T)^{1/2}}{(\log_2T)^{n+1/2}}.\]
Since
\[S_{n+1}(t+h)-S_{n+1}(t) = \int_t^{t+h} S_n(\tau)\mmd \tau,\]
it follows that
\[\max_{t\leq \tau\leq t+h}\pm S_n(\tau) \ge \pm \frac{1}{h}\big(S_{n+1}(t+h)-S_{n+1}(t)\big).\]
Since \(\beta^\prime>\beta\), we have \((T/2)^{\beta^\prime} \ge T^\beta\) and \(T/2+h \le T\) for all sufficiently large \(T\). Hence, \([t,t+h]\subseteq[T^\beta,T]\), which completes the proof of Corollary \ref{cor2}.

\section*{AI Declaration}
The authors used ChatGPT solely for language revision in this work and assume full responsibility for the mathematical arguments, results, and final text.
	
	\bibliographystyle{siam}
    \bibliography{references}
\end{document}